\documentclass[10pt]{article}

\usepackage{a4wide}
\usepackage{amssymb}
\usepackage{amsfonts}
\usepackage{amsmath}
\input xy
\xyoption{arrow} \xyoption{matrix}

\date{}

\newtheorem{proposition}{Proposition}[section]
\newtheorem{theorem}[proposition]{Theorem}
\newtheorem{lemma}[proposition]{Lemma}

\newtheorem{corollary}[proposition]{Corollary}

\def\der{\partial }

\def\nFM0{{\nu }_{F,M_0}}
\def\nFN0{{\nu }_{F,N_0}}
\def\nGN0{{\nu }_{G,N_0}}

\def\N0{ {\bf N}_0 }

\def\ra{\rightarrow}

\def\Xpm{X^{\pm }}

\def\s{\sigma}
\def\Z{\mathbb{Z}}

\def\l1{{\lambda}_1}

\def\a{\alpha}
\def\a0{ {\alpha }_0}
\def\a1{ {\alpha }_1}

\def\l{\lambda}
\def\o{\omega}

\def\nFGM0{{\nu }_{F,G,M_0}}

\def\nFN0{{\nu}_{F,N_0}}

\def\sm{{\sigma}^m}

\def\sm1{{\sigma}^{-1}}

\def\smtp1{{\sigma}^{-t+1}}

\def\o{\omega }
\def\S1{S^{-1}}

\def\Xpm1{X^{\pm 1}_1}

\def\sPM1{{\sigma }^{\pm 1}}
\def\sMP1{{\sigma }^{\mp 1 }}

\def\d{\delta}

\def\di{{\rm d.ind}}

\def\L{\Lambda}

\def\CA{{\cal A}}

\def\CD{{\cal D}}

\def\Ytm1{Y^{t-1}}
\def\Yim1{Y^{i-1}}

\def\CF{{\cal F}}
\def\CG{{\cal G}}
\def\CH{{\cal H}}

\def\Aut{{\rm Aut}}

\def\Der{{\rm Der }}
\def\ad{{\rm ad }}
\def\dim{{\rm dim }}

\def\ker{ {\rm ker } }

\def\CJ{ {\cal J}}

\def\D{ \Delta }

\def\SL2Z{ {\rm SL}_2({\bf Z}) }

\def\Gp1{ G^{1 , 1 } }
\def\P11{ P^{-1 , 1 } }
\def\Pp1{ P^{1 , 1 } }

\def\Supp{{\rm Supp}}

\def\nCLsr{{}^\nu\kern-2pt {\cal L}^{\sigma , \rho  }}
\def\nP{{}^\nu \kern-2pt P}
\def\nL{{}^\nu\kern-2pt L}
\def\nLL{{}^\nu\kern-2pt \Lambda}
\def\nPsr{{}^\nu\kern-2pt P^{\sigma , \rho  }}
\def\nLsr{{}^\nu\kern-2pt L^{\sigma , \rho  }}
\def\nuCL{{}^\nu\kern-2pt  {\cal L}}
\def\nCLsr{{}^\nu\kern-2pt {\cal L}^{\sigma , \rho  }}
\def\nCL1m{{}^\nu\kern-2pt {\cal L}^{-1 , 1  }}
\def\x1nu{x^\frac{1}{\nu}}
\def\xm1nu{x^{-\frac{1}{\nu}}}

\def\ra{\rightarrow }

\def\CB{{\cal B}}

\def\CH{ {\cal H}}

\def\nAM0{{\nu }_{{\cal A},M_0}}
\def\nAN0{{\nu }_{{\cal A},N_0}}

\def\End{ {\rm End }}
\def\Der{ {\rm Der }}
\def\CJ{ {\cal J }}

\def\det{ {\rm det }}
\def\ad{ {\rm ad }}

\def\ga{\mathfrak{a}}

\def\SL{{\rm SL}}

\def\di!{\frac{\der^i}{i!}}
\def\dik!{\frac{\der^k_i}{k!}}

\def\id{{\rm id}}

\def\N{\mathbb{N}}

\def\0{\overline{0}}
\def\1{\overline{1}}

\def\Ln1{\L_{n,\overline{1}}}

\def\a1{a_{\overline{1}}}

\def\St{{\rm St}}
\def\S{\Sigma}

\def\grad{{\rm grad}}

\def\vn1{\overrightarrow{n-1}}

\def\Sh{{\rm Sh}}

\def\Inn{{\rm Inn}}

\def\mJ{\mathbb{J}}
\def\mI{\mathbb{I}}

\def\mF{\mathbb{F}}

\def\mT{\mathbb{T}}

\def\mG{\mathbb{G}}
\def\mE{\mathbb{E}}

\def\K1{{\rm K}_1}

\def\hmI1{\widehat{\mI_1}}
\def\tmI1{\widetilde{\mI_1}}
\def\tmJ1{\widetilde{\mJ_1}}
\def\hB1{\widehat{B_1}}
\def\hCB1{\widehat{\CB_1}}

\def\ggu{\mathfrak{u}}

\def\Fix{{\rm Fix}}

\def\UAut{{\rm UAut}}

\def\mJ{\mathbb{J}}

\begin{document}

\author{V. V. \  Bavula   
}

\title{The group of automorphisms of the Lie algebra of derivations of a polynomial algebra}

\maketitle

\begin{abstract}

We prove that the group of automorphisms of the Lie algebra $\Der_K (P_n)$ of derivations of a polynomial algebra $P_n=K[x_1, \ldots , x_n]$ over a field of characteristic zero is canonically isomorphic to the the group of automorphisms of the polynomial algebra $P_n$.

$\noindent $



\end{abstract}


\section{Introduction}

In this paper, module means a left module, $K$ is a
field of characteristic zero and  $K^*$ is its group of units, and the following notation is fixed:
\begin{itemize}
\item $P_n:= K[x_1, \ldots , x_n]=\bigoplus_{\alpha \in \N^n}
Kx^{\alpha}$ is a polynomial algebra over $K$ where
$x^{\alpha}:=x_1^{\alpha_1}\cdots x_n^{\alpha_n}$,
 \item $G_n:=\Aut_K(P_n)$ is the group of automorphisms of the polynomial algebra $P_n$,
 \item $\der_1:=\frac{\der}{\der x_1}, \ldots , \der_n:=\frac{\der}{\der
x_n}$ are the partial derivatives ($K$-linear derivations) of
$P_n$,
\item    $D_n:=\Der_K(P_n) =\bigoplus_{i=1}^nP_n\der_i$ is the Lie
algebra of $K$-derivations of $P_n$ where $[\der , \d ]:= \der \d -\d \der $,
\item  $\d_1:=\ad (\der_1), \ldots , \d_n:=\ad (\der_n)$ are the inner derivations of the Lie algebra $D_n$ determined by the elements $\der_1, \ldots , \der_n$ (where $\ad (a)(b):=[a,b]$),
\item $\mG_n:=\Aut_{{\rm Lie}}(D_n)$ is the group of automorphisms of the Lie algebra $D_n$,
 \item $\CD_n:=\bigoplus_{i=1}^n K\der_i$,
 \item $\CH_n :=\bigoplus_{i=1}^n KH_i$ where $H_1:=x_1\der_1, \ldots , H_n:=x_n\der_n$,
 \item  $A_n:= K \langle x_1, \ldots
, x_n , \der_1, \ldots , \der_n\rangle  =\bigoplus_{\alpha , \beta
\in \N^n} Kx^\alpha \der^\beta$  is  the $n$'th {\em Weyl
algebra},
\item  for each natural number $n\geq 2$, $\ggu_n :=
K\der_1+P_1\der_2+\cdots +P_{n-1}\der_n$ is the  {\em Lie algebra of triangular polynomial derivations} (it is a Lie subalgebra of the Lie algebra $D_n$) and $ \Aut_K(\ggu_n)$ is its group of automorphisms.
    \end{itemize}

The aim of the paper is to prove the following theorem.

\begin{theorem}\label{11Mar13}
$\mG_n = G_n$.
\end{theorem}

{\it Structure of the proof}. (i) $G_n\subseteq \mG_n$ via the group monomorphism (Lemma \ref{b11Mar13}.(3))
$$G_n\ra \mG_n, \;\;  \s \mapsto \s : \der \mapsto \s (\der ):=\s \der \s^{-1}.$$


(ii) Let $\s \in \mG_n$. Then $\der_1':=\s (\der_1), \ldots , \der_n':=\s (\der_n)$ are commuting, locally nilpotent derivations of the polynomial algebra $P_n$ (Lemma \ref{c13Mar13}.(1)).

$\noindent $

(iii) $\bigcap_{i=1}^n \ker_{P_n}(\der_i') = K$  (Lemma \ref{c13Mar13}.(2)).

$\noindent $

(iv)(crux) There exists a polynomial automorphism $\tau \in G_n$ such that $\tau \s \in \Fix_{\mG_n}(\der_1, \ldots , \der_n)$ (Corollary \ref{b13Mar13}).

$\noindent $

(v) $\Fix_{\mG_n}(\der_1, \ldots , \der_n)=\Sh_n$ (Proposition \ref{B11Mar13}.(3)) where
$$\Sh_n:=\{ s_\l \in G_n\, | \, s_\l (x_1)=x_1+\l_1, \ldots , s_\l (x_n) = x_n+\l_n\}$$ is the {\em shift group} of automorphisms of the polynomial algebra $P_n$ and $\l = (\l_1, \ldots , \l_n)\in K^n$.

$\noindent $

(vi) By (iv) and (v), $\s \in G_n$, i.e. $\mG_n = G_n$.   $\Box $

$\noindent $

{\bf An analogue of the Jacobian Conjecture is true for $D_n$}. The Jacobian Conjecture claims that {\em certain} monomorphisms of the polynomial algebra $P_n$ are isomorphisms: {\em Every algebra endomorphism $\s $ of the polynomial algebra $P_n$ such that $\CJ (\s ):= \det (\frac{\der \s (x_i)}{\der x_j})\in K^*$ is an automorphism.} The condition that $\CJ (\s )\in K^*$ implies  that the endomorphism $\s$ is a monomorphism.

$\noindent $

{\bf Conjecture}. {\em  Every homomorphism of the Lie algebra $D_n$ is an automorphism.}


\begin{theorem}\label{10Mar12}
\cite{Bav-Lie-Un-MON} Every monomorphism of the Lie algebra $\ggu_n$ is an automorphism.
\end{theorem}

{\it Remark}. Not every epimorphism of the Lie algebra $\ggu_n$ is an automorphism. Moreover, there are countably many distinct ideals $\{ I_{i\o^{n-1}}\, | \, i\geq 0\}$ such that $$I_0=\{0\}\subset I_{\o^{n-1}}\subset I_{2\o^{n-1}}\subset \cdots \subset I_{i\o^{n-1}}\subset \cdots$$ and the Lie algebras $ \ggu_n/I_{i\o^{n-1}}$ and $\ggu_n$ are isomorphic (Theorem 5.1.(1), \cite{Bav-Lie-Un-GEN}).

Theorems \ref{10Mar12} and Conjecture have bearing of the Jacobian Conjecture and  the Conjecture of Dixmier  \cite{Dix}
 for the Weyl algebra $A_n$ over a field of characteristic zero that claims: {\em every homomorphism of the Weyl algebra is an automorphism}. The Weyl algebra $A_n$ is a simple algebra, so every algebra endomorphism of $A_n$  is a monomorphism. This conjecture is open since 1968 for all $n\geq 1$. It is stably equivalent to the Jacobian  Conjecture for the polynomial algebras as was shown by  Tsuchimoto
\cite{Tsuchi05}, Belov-Kanel and Kontsevich \cite{Bel-Kon05JCDP},
(see also  \cite{JC-DP} for a short proof which is based on the author's new inversion formula for polynomial automorphisms \cite{Bav-inform}).

$\noindent $

{\bf An analogue of the Conjecture of Dixmier is true for the algebra $\mI_1:= K\langle x, \frac{d}{dx}, \int \rangle$ of polynomial
 integro-differential operators}.

  \begin{theorem}\label{11Oct10}
{\rm (Theorem 1.1, \cite{Bav-cdixintdif})} Each algebra endomorphism of $\mI_1$ is an automorphism.
\end{theorem}

In contrast to the Weyl algebra $A_1=K\langle x, \frac{d}{dx} \rangle$, the algebra of polynomial
 differential operators, the algebra $\mI_1$ is neither a left/right Noetherian algebra nor a simple algebra. The left localizations, $A_{1,\der}$ and $\mI_{1, \der}$, of the algebras  $A_1$ and $\mI_1$ at the powers of the element $\der= \frac{d}{dx}$ are isomorphic. For the simple  algebra $A_{1,\der} \simeq \mI_{1, \der}$, there are algebra endomorphisms that are not automorphisms \cite{Bav-cdixintdif}.

$\noindent $

{\bf The group of automorphisms of the Lie algebra $\ggu_n$}.
In \cite{Bav-Lie-Un-AUT}, the group of automorphisms $\Aut_K(\ggu_n)$ of the Lie algebra $\ggu_n$ of triangular polynomial derivations is found ($n\geq 2$), it is isomorphic to an iterated semi-direct  product (Theorem 5.3,  \cite{Bav-Lie-Un-AUT}),
 $$\mT^n\ltimes (\UAut_K(P_n)_n\rtimes( \mF_n' \times  \mE_n ))  $$
 where $\mT^n$ is an algebraic  $n$-dimensional torus,   $\UAut_K(P_n)_n$ is an explicit factor group of the group $\UAut_K(P_n)$ of unitriangular polynomial automorphisms, $\mF_n'$ and $\mE_n$ are explicit groups that are isomorphic respectively to the groups $\mI$ and $\mJ^{n-2}$ where
    $\mI := (1+t^2K[[t]], \cdot )\simeq K^{\N}$ and
  $\mJ := (tK[[t]], +)\simeq K^\N$. Comparing the groups $G_n$ and $\Aut_K(\ggu_n)$ we see that the group $(\UAut_K(P_n)_n$ of polynomial automorphisms is a {\em tiny} part of the group $\Aut_K(\ggu_n)$ but in contrast $\mG_n = \Aut_K(P_n)$.
   It is shown that the {\em adjoint group} of automorphisms $\CA (\ggu_n)$  of the Lie algebra $\ggu_n $ is equal to the group $\UAut_K(P_n)_n$ (Theorem 7.1,  \cite{Bav-Lie-Un-AUT}). Recall that the {\em adjoint group} $\CA (\CG )$ of a Lie algebra $\CG$ is generated by the elements $ e^{ \ad (g)}:=\sum_{i\geq 0}\frac{\ad (g)^i}{i!}\in \Aut_K(\CG )$ where $g$ runs through all the locally nilpotent elements of the Lie algebra $\CG$ (an element $g$ is a {\em locally nilpotent element} if the inner derivation $\ad (g):= [g, \cdot ]$ of the Lie algebra $\CG$ is a locally nilpotent derivation).



\section{Proof of Theorem \ref{11Mar13} }\label{PPPAAA}

This section can be seen as a proof of Theorem \ref{11Mar13}. The proof is split into several statements that reflect `Structure of the proof of Theorem \ref{11Mar13}' given in the Introduction.

{\bf The Lie algebra $D_n$ is $\Z^n$-graded}. The Lie algebra
\begin{equation}\label{xadbd}
D_n =\bigoplus_{\alpha\in \N^n} \bigoplus_{i=1}^n Kx^\alpha \der_i
\end{equation}

 is a $\Z^n$-graded Lie algebra
$$D_n = \bigoplus_{\beta\in \Z^n} D_{n , \beta}\;\; {\rm where}\;\; D_{n,\beta}=\bigoplus_{\alpha -e_i=\beta}Kx^\alpha \der_i,$$
i.e. $[D_{n,\alpha} , D_{n,\beta }] \subseteq D_{n, \alpha +\beta}$ for all $\alpha ,\beta \in \N^n$ where $e_1:=(1, 0 , \ldots , 0), \ldots , e_n:=(0, \ldots , 0 , 1)$ is the canonical free basis for the free abelian group $\Z^n$. This follows from the commutation relations

\begin{equation}\label{xadbd1}
[x^\alpha\der_i, x^\beta \der_j]= \beta_i x^{\alpha+\beta - e_i} \der_j-\alpha_j x^{\alpha + \beta - e_j} \der_i.
\end{equation}
Clearly, for all $i,j=1, \ldots , n$ and $\alpha \in \N^n$,
\begin{equation}\label{xadbd2}
[H_j, x^\alpha \der_i]=\begin{cases}
\alpha_j x^{\alpha} \der_i & \text{if }j\neq i ,\\
(\alpha_i-1)x^{\alpha} \der_i& \text{if }j=i, \\
\end{cases}
\end{equation}
\begin{equation}\label{xadbd3}
[\der_j, x^\alpha \der_i]=\alpha_j x^{\alpha -e_j} \der_i.
\end{equation}
The {\em support} $\Supp (D_n):=\{ \beta \in \Z^n\, | \, D_{n,\beta}\neq 0\}$ is a submonoid of $\Z^n$. Let us find the support $\Supp (D_n)$, the graded components $D_{n,\beta}$ and their dimensions $\dim_K\, D_{n,\beta}$. For each $i=1, \ldots , n$, let $\N^{n,i}:=\{ \alpha \in \N^n \, | \, \alpha_i=0\}$ and $P_n^{\der_i}:=\ker_{P_n}(\der_i)$. It follows from the decompositions
$P_n = P_n^{\der_i}\oplus P_nx_i$ for $i=1, \ldots , n$ that
$$D_n = \bigoplus_{i=1}^n (P_n^{\der_i}\oplus P_nx_i)\der_i =\bigoplus_{i=1}^n P_n^{\der_i}\der_i \oplus \bigoplus_{i=1}^nP_nH_i, $$
\begin{equation}\label{Dnb2}
D_n =\bigoplus_{i=1}^n P_n^{\der_i}\der_i \oplus \bigoplus_{\alpha \in \N^n} x^\alpha \CH_n.
\end{equation}
Hence,
\begin{equation}\label{Dnb}
\Supp (D_n) =\coprod_{i=1}^n (\N^{n,i}-e_i) \coprod\N^n.
\end{equation}
\begin{equation}\label{Dnb1}
D_{n,\beta} =\begin{cases}
x^\alpha\der_i& \text{if }\beta = \alpha - e_i\in \N^{n,i}-e_i,\\
x^\beta \CH_n& \text{if }\beta \in \N^n.
\end{cases}
\end{equation}
$$\dim_K\, D_{n,\beta} =\begin{cases}
1& \text{if }\beta = \alpha - e_i\in \N^{n,i}-e_i,\\
n& \text{if }\beta \in \N^n.
\end{cases}$$

Let $\CG$ be a Lie algebra and $\CH$ be its Lie subalgebra. The {\em centralizer} $C_\CG (\CH ) := \{ x\in \CG \, | \, [ x, \CH ] =0\}$ of $\CH$ in $\CG$ is a Lie subalgebra of $\CG$. In particular, $Z(\CG ) := C_{\CG }(\CG ) $ is the {\em centre} of the Lie algebra $\CG$. The {\em normalizer} $N_\CG (\CH ) :=\{ x\in \CG \, | \, [ x, \CH ] \subseteq \CH\}$ of $\CH$ in $\CG$ is a Lie subalgebra of $\CG$, it is the largest Lie subalgebra of $\CG$ that contains $\CH $ as an ideal.

Let $V$ be a vector space over $K$. A $K$-linear map $\d : V\ra V$
is called a {\em locally nilpotent map} if $V=\bigcup_{i\geq 1} \ker
(\d^i)$ or, equivalently, for every $v\in V$, $\d^i (v) =0$ for
all $i\gg 1$. When  $\d$ is a locally nilpotent map in $V$ we
also say that $\d$ {\em acts locally nilpotently} on $V$. Every {\em nilpotent} linear map  $\d$, that is $\d^n=0$ for some $n\geq 1$, is a locally nilpotent map but not vice versa, in general.   Let
$\CG$ be a Lie algebra. Each element $a\in \CG$ determines the
derivation  of the Lie algebra $\CG$ by the rule $\ad (a) : \CG
\ra \CG$, $b\mapsto [a,b]$, which is called the {\em inner
derivation} associated with $a$. The set $\Inn (\CG )$ of all the
inner derivations of the Lie algebra $\CG$ is a Lie subalgebra of
the Lie algebra $(\End_K(\CG ), [\cdot , \cdot ])$ where $[f,g]:=
fg-gf$. There is the short exact sequence of Lie algebras
$$ 0\ra Z(\CG ) \ra \CG\stackrel{\ad}{\ra} \Inn (\CG )\ra 0,$$
that is $\Inn (\CG ) \simeq \CG / Z(\CG )$ where $Z(\CG )$ is the {\em centre} of the Lie algebra $\CG$ and $\ad ([a,b]) = [
\ad (a) , \ad (b)]$ for all elements $a, b \in \CG$. An element $a\in \CG$ is called a {\em locally nilpotent element} (respectively, a {\em nilpotent element}) if so is the inner derivation $\ad (a)$ of the Lie algebra $\CG$.

$\noindent $

{\bf The Cartan subalgebra $\CH_n$ of $D_n$}. A nilpotent Lie subalgebra $C$ of a Lie algebra $\CG$ is called a {\em Cartan subalgebra} of $\CG$ if it coincides with its normalizer. We use often the following obvious observation: {\em An abelian Lie subalgebra that coincides with its centralizer is a maximal abelian Lie subalgebra}.

\begin{lemma}\label{a11Mar13}
\begin{enumerate}
\item $\CH_n$ is a Cartan subalgebra of $D_n$.
\item  $\CH_n=C_{D_n}(\CH_n)$ is a maximal abelian subalgebra of $D_n$.
\end{enumerate}
\end{lemma}

{\it Proof}. Statements 1 and 2  follows from (\ref{Dnb}) and (\ref{Dnb1}). $\Box $

$\noindent $

{\bf $P_n$ is a $D_n$-module}. The polynomial algebra $P_n$ is a (left) $D_n$-module: $D_n \times P_n\ra P_n$, $(\der, p)\mapsto \der *p$. In more detail, if $\der = \sum_{i=1}^n a_i\der_i$ where $a_i\in P_n$ then
$$\der * p = \sum_{i=1}^n a_i\frac{\der p}{\der x_i}.$$
The field $K$ is a $D_n$-submodule of  $P_n$ and
\begin{equation}\label{IkerdiK}
\bigcap_{i=1}^n \ker_{P_n}(\der_i)= K.
\end{equation}

\begin{lemma}\label{xa11Mar13}
The $D_n$-module $P_n/K$ is simple with $\End_{D_n}(P_n/K)=K\id  $ where $\id$ is the identity map.
\end{lemma}

{\it Proof}. Let $M$ be a nonzero submodule of $P_n/K$ and $0\neq p\in M$. Using the actions of $\der_1, \ldots , \der_n$ on $p$ we obtain an element of $M$ of the form $\l x_i$ for some $\l\in K^*$. Hence, $x_i\in M$ and $x^\alpha = x^\alpha \der_i *x_i \in M$ for all $0\neq \alpha \in \N^n$. Therefore, $M = P_n /K$. Let $f\in \End_{D_n}(P_n/K)$. Then applying $f$ to the equalities $\der_i*(x_1+K)=\d_{i1}$ for $i=1, \ldots , n$, we obtain the equalities
$$ \der_i*f(x_1+K)=\d_{i1} \;\; {\rm for }\;\; i=1, \ldots , n.$$ Hence, $f(x_1+K)\in \bigcap_{i=2}^n \ker_{P_n/K}(\der_i) \cap \ker_{P_n/K}(\der_i^2) = (K[x_1]/K)\cap \ker_{P_n/K}(\der_i^2) =K(x_1+K)$. So, $f(x_1+K) = \l ( x_1+K)$ and so $f=\l \id$, by the simplicity of the $D_n$-module  $P_n/K$.

 $\Box $

$\noindent $


{\bf  The $G_n$-module $D_n$}. The Lie algebra $D_n$ is a $G_n$-module,
$$ G_n\times D_n\ra D_n, \;\; (\s , \der ) \mapsto \s (\der ) := \s \der \s^{-1}.$$
Every automorphism $\s \in G_n$ is uniquely determined by the elements
$$x_1':=\s (x_1), \; \ldots , \; x_n':=\s (x_n).$$
Let $M_n(P_n)$ be the algebra of $n\times n$ matrices over  $P_n$. The matrix  $J(\s) := (J(\s )_{ij}) \in M_n(P_n)$, where $J(\s )_{ij} =\frac{\der x_j'}{\der x_i}$,   is called the {\em Jacobian matrix} of the automorphism (endomorphism)  $\s$ and its determinant $\CJ (\s ) :=\det \, J(\s)$ is called the {\em Jacobian} of $\s$. So, the $j$'th column of $J(\s )$ is the {\em gradient} $\grad \, x_j':=(\frac{\der x_j'}{\der x_1}, \ldots , \frac{\der x_j'}{\der x_n})^T$  of the polynomial $x_j'$. Then the derivations
$$\der_1':= \s \der_1\s^{-1}, \; \ldots , \; \der_n':= \s\der_n\s^{-1}$$ are the partial derivatives of $P_n$ with respect to the variables $x_1', \ldots , x_n'$,
\begin{equation}\label{ddp=dxi}
\der_1'=\frac{\der}{\der x_1'}, \; \ldots , \; \der_n'=\frac{\der}{\der x_n'}.
\end{equation}
Every derivation $\der \in D_n$ is a unique sum $\der = \sum_{i=1}^n a_i\der_i$ where $a_i = \der *x_i\in P_n$. Let  $\der := (\der_1, \ldots , \der_n)^T$ and $ \der' := (\der_1', \ldots , \der_n')^T$ where $T$ stands for the transposition. Then
\begin{equation}\label{dp=Jnd}
\der'=J(\s )^{-1}\der , \;\; {\rm i.e.}\;\; \der_i'=\sum_{j=1}^n (J(\s )^{-1})_{ij} \der_j\;\; {\rm for }\;\; i=1, \ldots , n.
\end{equation}
In more detail, if $\der'=A\der $ where $A= (a_{ij})\in M_n(P_n)$, i.e. $\der_i=\sum_{j=1}^n a_{ij}\der_j$. Then for all $i,j=1, \ldots , n$,
$$\d_{ij}= \der_i'*x_j'=\sum_{k=1}^na_{ik}\frac{\der x_j'}{\der x_k}$$
where $\d_{ij}$ is the Kronecker delta function. The equalities above can be written in the matrix form as  $AJ(\s) = 1$ where $1$ is the identity matrix. Therefore, $A= J(\s )^{-1}$.

Suppose that a group $G$ acts on a set $S$. For a nonempty subset $T$ of $S$, $\St_G(T):=\{ g\in G\, | \, gT=T\}$ is the {\em stabilizer} of the set $T$ in $G$ and $\Fix_G(T):=\{ g\in G\, | \, gt=t$ for all $t\in T\}$  is the {\em fixator} of the set $T$ in $G$. Clearly, $\Fix_G(T)$ is a {\em normal} subgroup of $\St_G(T)$.

$\noindent $

{\bf The maximal abelian Lie subalgebra $\CD_n$ of $D_n$}.
\begin{lemma}\label{b11Mar13}

\begin{enumerate}
\item $C_{D_n}(\CD_n) =\CD_n$ and so $\CD_n$ is a maximal abelian Lie subalgebra of $D_n$.
\item $\Fix_{G_n}(\CD_n) = \Fix_{G_n}(\der_1, \ldots , \der_n) = \Sh_n$.
\item $D_n$ is a faithful $G_n$-module, i.e. the group homomorphism $G_n\ra \mG_n$, $ \s\mapsto \s : \der\mapsto \s \der \s^{-1}$,  is a monomorphism.
    \item $\Fix_{G_n}(\der_1, \ldots , \der_n , H_1, \ldots , H_n)=\{ e\}$.
\end{enumerate}
\end{lemma}

{\it Proof}. 1. Statement 1 follows from (\ref{xadbd1}).

2. Let $\s \in \Fix_{G_n}(D_n)$ and  $J(\s ) = (J_{ij})$. By (\ref{dp=Jnd}), $\der = J(\s ) \der $, and so, for all $i,j=1, \ldots , n$,
$\d_{ij} = \der_i*x_j=J_{ij}$,
i.e. $J(\s ) = 1$, or equivalently, by (\ref{IkerdiK}),
$$x_1'=x_1+\l_1, \ldots , x_n'=x_n+\l_n$$
for some scalars $\l_i\in K$, and so $\s\in \Sh_n$.

3 and 4. Let $\s\in \Fix_{G_n}=(\der_1, \ldots , \der_n , H_1, \ldots , H_n)$. Then $\s \in  \Fix_{G_n}(\der_1, \ldots , \der_n)=\Sh_n$, by statement 2. So, $\s (x_1) = x_1+\l_1, \ldots , \s (x_n) = x_n+\l_n$ where $\l_i\in K$. Then $x_i\der_i = \s (x_i\der_i) = (x_i+\l_i) \der_i$ for $i=1, \ldots , n$, and so $\l_1=\cdots = \l_n=0$. This means that $\s = e$. So,
 $\Fix_{G_n}=(\der_1, \ldots , \der_n , H_1, \ldots , H_n)=\{ e\}$ and $D_n$ is a faithful $G_n$-module.  $\Box $

$\noindent $

By Lemma \ref{b11Mar13}.(3), we identify the group $G_n$ with its image in $\mG_n$.

\begin{lemma}\label{c11Mar13}

\begin{enumerate}
\item $D_n$ is a simple Lie algebra.
\item $Z(D_n)=\{ 0\}$.
\item $[D_n, D_n]=D_n$.
\end{enumerate}
\end{lemma}

{\it Proof}. 1. Let $0\neq a\in D_n$ and $\ga = (a)$ be the ideal of the Lie algebra $D_n$ generated by the element $a$. We have to show that $\ga = D_n$. Using the inner derivations $\d_1, \ldots , \d_n$  we see that $\der_i\in \ga$ for some $i$. Then $\ga = D_n$ since
$$ x^\alpha \der_j = (\alpha_i+1)^{-1} [ \der_i, x^{\alpha +e_i}\der_j]\in \ga$$ for all $\alpha$ and $j$.

2 and 3. Statements 2 and 3 follow from statement 1.
 $\Box $

$\noindent $

\begin{proposition}\label{B11Mar13}

\begin{enumerate}
\item $\Fix_{\mG_n} (\der_1, \ldots , \der_n, H_1, \ldots , H_n) = \{ e\}$.
\item Let $\s, \tau \in \mG_n$. Then $\s = \tau $ iff $\s (\der_i) = \tau (\der_i)$ and $\s (H_i) = \tau (H_i)$ for $i=1, \ldots , n$.
\item $\Fix_{\mG_n} (\der_1, \ldots , \der_n) = \Sh_n$.
\end{enumerate}
\end{proposition}

{\it Proof}. 1. Let $\s\in F:= \Fix_{\mG_n} (\der_1, \ldots , \der_n, H_1, \ldots , H_n)$. We have to show that $\s = e$. Since $\s \in \Fix_{\mG_n} (H_1, \ldots , H_n)$, the automorphism $\s$ respects  the weight decomposition  of $D_n$. By (\ref{Dnb1}),
 $\s (x^\alpha \der_i) = \l_{\alpha , i} x^\alpha \der_i$ for all $\alpha \in \N^{n,i}$ and $i=1, \ldots , n$ where $\l_{\alpha , i}\in K$.
 Clearly, $\l_{0 , i}=1$ for $i=1, \ldots , n$. Since $\s \in \Fix_{\mG_n} (\der_1, \ldots , \der_n)$, by applying $\s$ to the relations $\alpha_jx^{\alpha-e_j}\der_i=[\der_j, x^\alpha\der_i]$, we get the relations
 $$\alpha_j \l_{\alpha-e_j,i} x^{\alpha-e_j} \der_i= [ \der_j, \l_{\alpha, i} x^\alpha\der_i]=  \alpha_j\l_{\alpha , i}x^{\alpha - e_j}\der_i.$$
Hence $\l_{\alpha , i} = \l_{\alpha - e_j, i}$ provided $\alpha_j\neq 0$. We conclude that all the coefficients $\l_{\alpha , i}$ are equal to one of the coefficients $\l_{e_i, j}$ where $i,j=1, \ldots , n$ and $i\neq j$. The relations $\der_j=[\der_i, x_i\der_j]$ implies the relations  $\der_j= [\der_i, \l_{e_i, j}x_i\der_j]=\l_{e_i, j}\der_j$, hence all the coefficients $\l_{e_i, j}$  are equal to 1.
 So, $\oplus_{i=1}^n P_n^{\der_i} \der_i\subseteq \CF := \Fix_{D_n}(\s ) :=\{ \der\in D_n \,| \, \s (\der )  = \der\}$. To finish the proof of statement 1 it suffices to show that $x^\alpha H_i\in \CF$ for all $\alpha \in \N^n$ and $i=1, \ldots , n$, see  (\ref{Dnb2}) and (\ref{Dnb}). We use induction on $|\alpha |:=\alpha_1+\cdots + \alpha_n$. If $|\alpha |=0$ the statement is obvious as $\s\in F$.  Suppose that $|\alpha |>0$. Using the commutation relations
\begin{equation}\label{djxaHj}
[\der_j, x^\alpha H_i] = \begin{cases}
\alpha_jx^{\alpha - e_j}H_i& \text{if }j\neq i,\\
(\alpha_i+1) x^\alpha \der_i& \text{if }j=i, \\
\end{cases}
\end{equation}
 the induction and the previous case,  we see that
$$ [\der_j, \s (x^\alpha H_i) - x^\alpha H_i]=0 \;\; {\rm for }\;\; i=1, \ldots , n. $$
Therefore, $  \s (x^\alpha H_i) - x^\alpha H_i\in C_{D_n} (\CD_n) = \CD_n$. Since the automorphism $\s$ respects the weight decomposition of $D_n$, we must have  $  \s (x^\alpha H_i) - x^\alpha H_i\in x^\alpha \CH_n \cap \CD_n=\{ 0 \}$. Hence, $x^\alpha H_i\in \CF$, as required.

2. Statement 2 follows from statement 1.

3. Clearly, $\Sh_n \subseteq F=\Fix_{\mG_n} (\der_1, \ldots , \der_n)$. Let $\s \in F$ and $H_i':=\s (H_i), \ldots , H_n':=\s (H_n)$. Applying the automorphism $\s$ to the commutation relations $[\der_i, H_j]=\d_{ij}\der_i$ gives the relations $[\der_i, H_j']=\d_{ij}\der_i$. By taking the difference, we see that $[\der_i, H_j'-H_j]=0$ for all $i$ and $j$. Therefore, $H_i'=H_i+d_i$ for some elements  $d_i\in C_{D_n}(\CD_n) = \CD_n$ (Lemma \ref{b11Mar13}.(1)),  and so $d_i=\sum_{j=1}^n \l_{ij}\der_j$ for some elements $\l_{ij}\in K$. The elements $H_1', \ldots , H_n'$ commute, hence
$$ [ H_j, \der_i]= [H_i, \der_j] \;\; {\rm for \; all}\;\; i,j,  $$
or equivalently,
$$ \l_{ij}\der_j= \l_{ji}\der_i \;\; {\rm for \; all}\;\; i,j.  $$
This means that $\l_{ij}=0$ for all $i\neq j$, i.e. $$H_i'= H_i+\l_{ii}\der_i= (x_i+\l_{ii})\der_i= s_\l (H_i)$$
where $s_\l \in \Sh_n$, $s_\l (x_i) = x_i+\l_{ii}$ for all $i$. Then $s_\l^{-1}\s \in \Fix_{\mG_n} (\der_1, \ldots , \der_n, H_1, \ldots , H_n) = \{ e\}$  (statement 2), and so $\s = s_\l \in \Sh_n$. $\Box $






\begin{lemma}\label{c13Mar13}
Let $\s \in \mG_n$ and $\der_1':=\s (\der_1), \ldots , \der_n':= \s (\der_n)$. Then
\begin{enumerate}
\item $\der_1', \ldots , \der_n'$ are commuting, locally nilpotent derivations of $P_n$.
\item $\bigcap_{i=1}^n\ker_{D_n}(\der_i')=K$.
\end{enumerate}
\end{lemma}

{\it Proof}. 1. The derivations $\der_1', \ldots , \der_n'$  commute  since $\der_1, \ldots , \der_n$ are commute. The inner derivations $\d_1, \ldots , \d_n$  of the Lie algebra $D_n$ are commuting and locally nilpotent.  Hence, inner derivations
 $$\d_1':= \ad (\der_1'), \ldots , \d_n':= \ad (\der_n')$$
 of the Lie algebra $D_n$ are commuting and  locally nilpotent. The vector space $P_n\der_i'$ is closed under the derivations $\d_j'$ since
 $$ \d_j'(P_n\der_i') = [ \der_j', P_n\der_i']= (\der_j' *P_n)\cdot \der_i' \subseteq P_n\der_i'.$$
 Therefore, $\der_1', \ldots , \der_n'$ are locally nilpotent derivations of the polynomial algebra $P_n$.

 2.   Let $\l \in \bigcap_{i=1}^n \ker_{P_n}(\der_i')$. Then
 $$\l \der_1' \in C_{D_n}(\der_1', \ldots , \der_n')=\s (C_{D_n}(\der_1, \ldots , \der_n))= \s (C_{D_n}(\CD_n))= \s (\CD_n) = \s (\bigoplus_{i=1}^n K\der_i) = \bigoplus_{i=1}^n K\der_i', $$
since $C_{D_n}(\CD_n)=\CD_n$, Lemma \ref{b11Mar13}.(1). Then $\l \in K$ since otherwise the infinite dimensional space $\bigoplus_{i\geq 0} K\l^i \der_1'$ would be a subspace of a finite dimensional space $\s (\CD_n)$.  $\Box $

$\noindent $

The following lemma is well-known and it is easy to prove.

\begin{lemma}\label{Aslice}
Let $\der$ be a locally nilpotent derivation of a commutative $K$-algebra $A$ such that $\der (x) =1$ for some element $x\in A$. Then $A= A^\der [x]$ is a polynomial algebra over the ring $A^\der := \ker (\der )$ of constants of the derivation $\der$ in the variable $x$.
\end{lemma}

The next theorem is the most important point in the proof of Theorem \ref{11Mar13} and, roughly speaking,  the main reason why Theorem \ref{11Mar13} holds.
\begin{theorem}\label{A13Mar13}
Let $\der_1', \ldots , \der_n'$ be commuting, locally nilpotent derivations of the polynomial algebra $P_n$ such that $\bigcap_{i=1}^n \ker_{P_n}(\der_i')=K$. Then there exist polynomials $x_1', \ldots , x_n'\in P_n$ such that
\begin{equation}\label{con*}
\der_i'*x_j'=\d_{ij}.
\end{equation}
Moreover, the algebra homomorphism
$$\s : P_n\ra P_n , \;\; x_1\mapsto x_1', \ldots , x_n\mapsto x_n'$$ is an automorphism such that $\der_i'= \s \der_i \s^{-1} = \frac{\der}{\der x_i'}$ for $i=1, \ldots , n$.
\end{theorem}

{\it Proof}. Case $n=1$: By Lemma \ref{c13Mar13}, the derivation $\der_1'$ of the polynomial algebra $P_1$ is a locally nilpotent derivation with $K_1':=\ker_{P_1}(\der_1') = K$. Hence, $\der_1'*x_1'=1$ for some polynomial $x_1'\in P_1$. By Lemma \ref{Aslice}, $K[x_1]= K_1'[x_1']= K[x_1']$, and so $\s : K[x_1]\ra K[x_1]$, $ x\mapsto x_1'$, is an automorphism such that $ \der_1' = \frac{d}{d x_1'} = \s \frac{d}{d x_1}\s^{-1}$.

Case $n\geq 2$. Let $K_i':= \ker_{P_n}(\der_i')$ for $i=1, \ldots, n$. Clearly, $ K\subseteq K_i'$.

(i) $K_i'\neq K$ {\em for} $i=1, \ldots , n$: If $K_i'= K$ for some $i$ then by the same argument as in the case $n=1$ there exists a polynomial $x_i'\in P_n$ such that $\der_i'*x_i'=1$, and so $P_n= K_i'[x_i'] = K[x_i]$, a contradiction.

(ii) Let $m $ be the maximum  of ${\rm card} (I)$ such $\emptyset \neq I \subseteq \{ 1, \ldots , n-1\}$ and $\bigcap_{i\in I} K_i'\neq K$. By (i), $2\leq m \leq n-1$. Changing (if necessary) the order of the derivations $\der_1',\ldots , \der_n'$  we may assume that $A:= \bigcap_{i=1}^m K_i'\neq K$.  Then the algebra $A$ is infinite dimensional (since $K\neq A \subseteq P_n$) and invariant under the action of the derivations $\der_j'$ for $j=m+1, \ldots , n$.  By the choice of $m$,
$$ A^{\der_j'}= K_j'\cap \bigcap_{i=1}^m K_i' = K\;\; {\rm for}\;\; j=m+1, \ldots , n$$
and the derivations $\der_j'$ acts locally nilpotently on the algebra $A^{\der_j'}$. Therefore, for each index $j=m+1, \ldots , n$, there exists an element $x_j'\in A$ such that $\der_j'*x_j'=1$, and so (Lemma \ref{Aslice})
\begin{equation}\label{Adjr}
A= A^{\der_j'} [ x_j']= K[x_j']\;\; {\rm for}  \;\; j=m+1, \ldots , n.
\end{equation}
(ii)(a) Suppose that $m=n-1$, i.e. $\der_i'*x_n'=\d_{in}$ for all $i=1, \ldots , n$. By Lemma \ref{Aslice}, $P_n = K_n'[x_n']$. The algebra $K_n'$ admits the set of commuting,  locally nilpotent derivations
$$ \der_1'':= \der_1'|_{K_n'}, \ldots , \der_{n-1}'':= \der_{n-1}'|_{K_n'}$$
with $\bigcap_{i=1}^{n-1} \ker_{K_n'} (\der_i'') = K_n'\cap \bigcap_{i=1}^{n-1} K_i' = K.$

(ii)(b) Suppose that $m<n-1$. By (\ref{Adjr}),
$$ K^*x_{m+1}' +K= K^*x_{m+2}' +K=\cdots = K^*x_n ' +K, $$
and so $\l_j := \der_j'*x_n'\in K$ for $j=m+1, \ldots , n-1$. Hence, $(\der_j'-\l_j\der_n') *x_n'=0$ for $j=m+1, \ldots , n-1$. A linear combination of commuting, locally nilpotent derivations  is a locally nilpotent derivation (the proof boils down to the case $\der +\d$ of two commuting, locally nilpotent derivations, then the result follows from $(\der +\d )^m = \sum_{i=0}^m {m \choose i} \der^i \d^{m-i}$ and $ \der^i \d^{m-i}= \d^{m-i} \der^i$).
Using the set of commuting, locally nilpotent derivations $\der_1', \ldots , \der_n'$ that satisfy (\ref{con*}) we obtain  the set of commuting, locally nilpotent derivations
$$\d_1':=\der_1',\;  \ldots , \; \d_m':=\der_m',\;  \d_{m+1}':=\der_{m+1}'-\l_{m+1}\der_n', \; \ldots , \;  \d_{n-1}':=\der_{n-1}'-\l_{n-1}\der_n',\;  \d_n':=\der_n$$
that satisfy (\ref{con*}) with
$$ \d_i'*x_n'=\d_{in}\;\; {\rm for}\;\; i=1, \ldots , n.$$ Then repeating the arguments of (ii)(a), we see that $P_n = K_n'[x_n']$. The algebra $K_n'$ admits the set of commuting, locally nilpotent derivations
$$ \der_1'':= \d_1'|_{K_n'},\;  \ldots , \; \der_{n-1}'':= \d_{n-1}'|_{K_n'}$$
with $$\bigcap_{i=1}^{n-1} \ker_{K_n'} (\der_i'') = K_n'\cap \bigcap_{i=1}^{n-1}\ker_{P_n}(\d_i')=  K_n'\cap \bigcap_{i=1}^{n-1}\ker_{P_n}(\der_i')=
\bigcap_{i=1}^n K_i' = K.$$
(iii) Using the  cases (ii)(a) and (ii)(b) $n-1$ more times we find polynomials $x_1', \ldots , x_n'$ and commuting set of locally nilpotent derivations of $P_n$,  say,  $\D_1, \ldots , \D_n$ that satisfy (\ref{con*}) and such that

$(\alpha)$ $\D_i*x_j'=\d_{ij}$ for all $i,j=1, \ldots , n$;

$(\beta )$ the $n$-tuple of derivations $\D = (\D_1, \ldots , \D_n)^T$ is obtained from the
 $n$-tuple of derivations $\der'= (\der_1' , \ldots , \der_n')^T$ by  unitriangular (hence invertible) scalar  matrix $\L = (\l_{ij})\in M_n(K)$ such that $ \D = \L \der'$; and

 $(\gamma )$ (where $K_1'':=\ker_{P_n}(\D_1),  \ldots , K_n'':=\ker_{P_n}(\D_n)$)
 \begin{eqnarray*}
P_n &=& K_n''[x_n']= (K_{n-1}''\cap K_n'')[x_{n-1}', x_n']= \cdots = (\bigcap_{i=s}^n K_i'') [x_s',\ldots ,  x_n']= \cdots \\
&=& (\bigcap_{i=1}^n K_i'') [x_1',\ldots ,  x_n']
=K[x_1',\ldots ,  x_n'].
 \end{eqnarray*}
(iv) Replacing the row $x'= (x_1', \ldots , x_n')$ by the row $x'\L$ gives the required elements of the theorem. Indeed,  by $(\alpha )$, $\L \cdot (\der_i'*x_j') = 1$, the identity $n\times n$ matrix. Hence, $ (\der_i'*x_j') \cdot \L = 1$, as required.

(v) Let $x_1' , \ldots , x_n'$ be the set of polynomials as in the theorem. Then $\s $ is an algebra automorphism (see $(\gamma )$ and (iv)) such that $\der_i'= \s \der_i\s^{-1} = \frac{\der}{\der x_i'}$ for $i=1, \ldots , n$. $\Box$

\begin{corollary}\label{b13Mar13}
Let $\s \in \mG_n$. Then $\tau \s \in \Fix_{\mG_n} (\der_1, \ldots , \der_n)$ for some $\tau \in G_n$.
\end{corollary}

{\it Proof}. By Lemma \ref{c13Mar13}, the elements $\der_1':= \s (\der_1), \ldots , \der_n':= \s (\der_n)$ satisfy the assumptions of Theorem \ref{A13Mar13}. By Theorem \ref{A13Mar13},  $\der_1':= \tau^{-1}  (\der_1), \ldots , \der_n':= \tau^{-1}  (\der_n)$ for some $\tau \in G_n$. Therefore, $ \tau \s \in \Fix_{\mG_n} (\der_1, \ldots , \der_n)$. $\Box $

$\noindent $

{\bf Proof of Theorem \ref{11Mar13}}. Let $\s \in \mG_n$. By Corollary \ref{b13Mar13}, $\tau \s \in \Fix_{\mG_n} (\der_1, \ldots , \der_n)=\Sh_n$ (Proposition \ref{B11Mar13}.(3)). Therefore, $\s \in G_n$, i.e. $\mG_n = G_n$. $\Box$


$${\bf Acknowledgements}$$

 The work is partly supported by  the Royal Society  and EPSRC.

\small{

Department of Pure Mathematics

University of Sheffield

Hicks Building

Sheffield S3 7RH

UK

email: v.bavula@sheffield.ac.uk}

\end{document}